\documentclass[11pt]{article}
\usepackage{amssymb}
\usepackage{mathrsfs}
\usepackage{tikz}
\usepackage{mathrsfs}
\usepackage{amssymb}
\usepackage{amsmath}
\usepackage{bm}
\usepackage{hyperref}
\usepackage{graphicx}
\usepackage{pst-poly}  
\usepackage{pst-plot}  

\newtheorem{thm}{Theorem}[section]
\newtheorem{lem}[thm]{Lemma}
\newtheorem{Def}[thm]{Definition}
\newtheorem{cor}[thm]{Corollary}

\newtheorem{rem}[thm]{Remark}

\newtheorem{prob}[thm]{Problem}

\newenvironment{pf}[1][Proof]{\noindent\textbf{#1.} }{\hfill\rule{2mm}{2mm}}

\makeatletter \@addtoreset{equation}{section} \makeatother

\begin{document}

\title
{Connectivity and eigenvalues of graphs with given girth \\
or clique number}

\author
{Zhen-Mu Hong$\,{\rm ^a}$, Hong-Jian Lai$\,{\rm ^b}$\thanks{E-mail addresses: zmhong@mail.ustc.edu.cn
(Zhen-Mu Hong), hjlai@math.wvu.edu (Hong-Jian Lai), xzj@mail.ustc.edu.cn (Zheng-Jiang Xia).}, Zheng-Jiang Xia$\,{\rm ^a}$\\
\\
{\footnotesize ${\rm ^a}$School of Finance, Anhui University of Finance \& Economics, Bengbu, Anhui 233030, China} \\
{\footnotesize ${\rm ^b}$Department of Mathematics, West Virginia University, Morgantown, WV 26506, USA}}


\date{}
\maketitle


\begin{center} {\bf Abstract} \end{center}

Let $\kappa'(G)$, $\kappa(G)$, $\mu_{n-1}(G)$ and $\mu_1(G)$ denote the edge-connectivity,
vertex-connectivity, the algebraic connectivity and the Laplacian spectral radius of $G$, respectively.
In this paper, we prove that for integers $k\geq 2$ and $r\geq 2$, and any simple graph $G$ of order $n$ with
minimum degree $\delta\geq k$, girth $g\geq 3$ and clique number $\omega(G)\leq r$,
the edge-connectivity $\kappa'(G)\geq k$ if $\mu_{n-1}(G) \geq  \frac{(k-1)n}{N(\delta,g)(n-N(\delta,g))}$
or if $\mu_{n-1}(G) \geq  \frac{(k-1)n}{\varphi(\delta,r)(n-\varphi(\delta,r))}$,
where $N(\delta,g)$ is the Moore bound on the smallest possible number of vertices
such that there exists a $\delta$-regular simple graph with girth $g$,
and $\varphi(\delta,r) = \max\{\delta+1,\lfloor\frac{r\delta}{r-1}\rfloor\}$.
Analogue results involving $\mu_{n-1}(G)$ and $\frac{\mu_1(G)}{\mu_{n-1}(G)}$ to
characterize vertex-connectivity of graphs with fixed girth and clique number are also presented.
Former results in [Linear Algebra Appl. 439 (2013) 3777--3784], [Linear Algebra Appl. 578 (2019) 411--424],
[Linear Algebra Appl. 579 (2019) 72--88], [Appl. Math. Comput. 344-345 (2019) 141--149] and
[Electronic J. Linear Algebra 34 (2018) 428--443] are improved or extended.

\vskip10pt




\section{Introduction}


We only consider finite and simple graphs in this paper.
Undefined notation and terminologies will follow Bondy and Murty \cite{BoMu08}.
Let $G=(V,E)$ be a graph of order $n$. We use $\kappa(G)$, $\kappa'(G)$,
$\delta(G)$ and $\Delta(G)$ to denote the vertex-connectivity, the edge-connectivity,
the minimum degree and the maximum degree of a graph $G$, respectively.
The girth $g(G)$ of a graph $G$ is the length of a shortest cycle in $G$ if
it contains at least one cycle, and $g(G)=\infty$ if $G$ is acyclic.
A clique of a graph is a set of mutually adjacent vertices, and that the maximum
size of a clique of a graph $G$, the clique number of $G$, is denoted $\omega(G)$.
For a vertex subset $S\subseteq V(G)$, $G[S]$ is the subgraph of $G$ induced by $S$.

Let $G=(V,E)$ be a simple graph with vertex set
$V=V(G)=\{v_1,v_2,\dots,v_n\}$ and edge set $E=E(G)$. The
adjacency matrix of $G$ is defined to be a $(0,1)$-matrix
$A(G)=(a_{ij})_{n \times n}$, where $a_{ij}=1$ if $v_i$ and $v_j$ are adjacent,
$a_{ij}=0$ otherwise. As $G$ is simple and undirected, $A(G)$ is a
symmetric $(0, 1)$-matrix. The adjacency eigenvalues of $G$ are the eigenvalues of $A(G)$.
Denoted by $D(G) = {\rm diag}\{d_G(v_1), d_G(v_2), \dots , d_G(v_n)\}$,
the degree diagonal matrix of $G$, where $d_G(v_i)$ denotes the degree of $v_i$.
The matrices $L(G) = D(G) - A(G)$ and $Q(G) = D(G) + A(G)$ are called
the Laplacian matrix and the signless Laplacian matrix of $G$, respectively.
We use $\lambda_i(G)$, $\mu_i(G)$ and $q_i(G)$ to denote the $i$th largest
eigenvalue of $A(G)$, $L(G)$ and $Q(G)$, respectively.

The second smallest Laplacian eigenvalue $\mu_{n-1}(G)$ is called algebraic connectivity by Fiedler \cite{Fied73,Fied75}.
Fiedler \cite{Fied73} initiated the investigation on the relationship between graph connectivity
and graph eigenvalues, and showed that $\mu_{n-1}(G)\leq \kappa(G)\leq \kappa'(G)$.
Kirkland et al. \cite{KMNS02} investigated the graphs with equal algebraic connectivity and vertex-connectivity.
It is worth to mention that Cioab\u{a} in \cite{Cioa10} investigated the relationship
between edge-connectivity and adjacency eigenvalues of regular graphs.
From then on, the edge-connectivity problem has been intensively studied by many researchers, as found in
\cite{ABMO18,Cioa10,CiWo12,GuXi13,GLLY16,LiSh13,LiHL14,LHGL14,LiLT19,LiLT13,OSui16}, among others.
For the vertex-connectivity of graphs, one can refer to \cite{ABMO18,HoXL19,LLTW19,OSar16}.
In \cite{ABMO18}, Abiad et al. raised the following research problem.

\begin{prob}{(Abiad et al. \cite{ABMO18})}\label{problem1.1}
For a $d$-regular simple graph or multigraph $G$ and for $2\leq k \leq d$,
what is the best upper bound on $\lambda_2(G)$ which guarantees $\kappa'(G)\geq k$ or $\kappa(G)\geq k$ ?
\end{prob}

A number of results are related to Problem \ref{problem1.1}, as shown in the following theorem.

\begin{thm}\label{thm-recent}
Let $d,k$ be integers with $d\geq k\geq 2$, and let $G$ be a simple graph of order $n$
with minimum degree $\delta\geq k$.

(i) (Cioab\v{a} \cite{Cioa10}) If $G$ is $d$-regular and
$\lambda_2(G) \leq d - \frac{(k-1)n}{(d+1)(n-d-1)}$, then $\kappa'(G)\geq k$.

(ii) (Li and Shi \cite{LiSh13}, Liu et al. \cite{LiHL14})
If $\lambda_2(G) \leq \delta - \frac{(k-1)n}{(\delta+1)(n-\delta-1)}$, then $\kappa'(G)\geq k$.

(iii) (Liu et al. \cite{LiLT13}) If $\mu_{n-1}(G) \geq \frac{(k-1)n}{(\delta+1)(n-\delta-1)}$ or
$q_2(G) \leq 2\delta-\frac{(k-1)n}{(\delta+1)(n-\delta-1)}$, then $\kappa'(G)\geq k$.

(iv) (Abiad et al. \cite{ABMO18}) Let $G$ be a $d$-regular graph.
If $k\geq 3$ and $\lambda_2(G) < d - \frac{(k-1)d n}{2(d-k+2)(n-d+k-2)}$,
then $\kappa(G)\geq k$. If $\lambda_2(G) < d - \frac{ d n}{2(d+1)(n-d-1)}$,
then $\kappa(G)\geq 2$.
\end{thm}

As can be seen in \cite{HoXL19} or will be seen in Section 4, for any real number $p>0$,
if $q_2(G)\leq 2\delta(G)-p$ or $\lambda_2(G)\leq \delta(G)-p$, then $\mu_{n-1}(G)\geq p$.
Moreover, it is known that if $\mu_{n-1}(G)>0$, then $\kappa'(G)\geq \kappa(G)\geq 1$.
Therefore, we focus on establishing the lower bounds on $\mu_{n-1}(G)$ which guarantee $\kappa'(G)\geq k$
or $\kappa(G)\geq k$. By Theorem \ref{thm-recent}, it is natural to discuss
Problem \ref{problem1.1} for bipartite graphs or triangle-free graphs and drop the graph regularity.
Note that triangle-free graphs have girth at least 4, or equivalently clique number at most 2. Thus,
to get better lower bounds on algebraic connectivity, we consider
graphs with fixed girth or clique number. In this paper, we improve or extend some recent results.
In order to state some known results, we need the the following definition.


\begin{Def}\label{def1.1}
For integers $\delta, g$ with $\delta\geq 2$ and $g\geq 3$, let $t=\lfloor\frac{g-1}{2}\rfloor$. Define
$$
N(\delta,g) :=
\left\{
\begin{array}{ll}
  1+\delta\sum_{i=0}^{t-1}(\delta-1)^i, &  \text{if~} g=2t+1; \\
  2\sum_{i=0}^{t}(\delta-1)^i, &  \text{if~} g=2t+2.
\end{array}
\right.
$$
\end{Def}

Tutte \cite{Tutt47} initiated the cage problem, which seeks, for any given integers $d$ and $g$ with $d\geq 2$ and
$g\geq 3$, the smallest possible number of vertices $n(d, g)$ such that there exists a $d$-regular simple
graph with girth $g$. $N(d,g)$ in Definition \ref{def1.1} is a tight lower bound
(often called the Moore bound) on $n(d, g)$ which can be found in \cite{ExJa11}.

The results in Theorem \ref{thm-recent} have been improved or extended in \cite{LiLT13,LiLT19,LLTW19,HoXL19} as follows.

\begin{thm}(Liu et al. \cite{LiLT13})\label{thm-LiLT13}
Let $k\geq 2$ be an integer, and $G$ be a connected graph of order $n$ with girth $g\geq 3$ and minimum degree
$\delta\geq k$. If $\mu_{n-1}(G)\geq \frac{(k-1)n}{g(n-g)}$, then $\kappa'(G)\geq k$.
Moreover, if $\delta\geq 3$ and $\mu_{n-1}(G)\geq \frac{(k-1)n}{\frac{4}{9}N(\delta,g)(n-\frac{4}{9}N(\delta,g))}$,
then $\kappa'(G)\geq k$.
\end{thm}

\begin{thm}(Liu et al. \cite{LiLT19})\label{thm-LiLT19}
Let $k\geq 2$ be an integer, and $G$ be a connected graph of order $n$ with girth $g\geq 3$ and minimum degree
$\delta\geq k$. Let $f(2,g)=g$, $t=\lfloor\frac{g-1}{2}\rfloor$ and for $\delta\geq 3$
$
f(\delta,g)= N(\delta,g)-\sum\nolimits_{i=1}^{t-1}(\delta-1)^i.
$
If $\mu_{n-1}(G)\geq \frac{(k-1)n}{f(\delta,g)(n-f(\delta,g))}$,
then $\kappa'(G)\geq k$.
\end{thm}

\begin{thm}(Liu et al. \cite{LLTW19})\label{thm-LLTW19}
Let $k\geq 2$ be an integer, and $G$ be a connected graph of order $n$ with
maximum degree $\Delta$, minimum degree $\delta\geq k$, girth
$g\geq 3$. Let $t = \lfloor\frac{g-1}{2}\rfloor$ and
\begin{equation*}
\nu(\delta,g,k)
=\left\{
\begin{array}{ll}
N(\delta,g)-(k-1)\sum\nolimits_{i=0}^{t-1}(\delta-1)^i, & \mbox{if $g=2t+1$, or $g=2t+2$ and $\delta\geq 3$}; \\
2t+1, & \mbox{if $g=2t+2$ and $\delta=2$}.
\end{array}
\right.
\end{equation*}
If $\mu_{n-1}(G)\geq \frac{(k-1)n\Delta}{2\nu(\delta,g,k)(n-\nu(\delta,g,k))}$, then $\kappa(G)\geq k$.
\end{thm}


\begin{thm}(Hong et al. \cite{HoXL19})\label{thm-HoXL19}
Let $k$ be an integer and $G$ be a simple graph of order $n$ with maximum degree $\Delta$
and minimum degree $\delta\geq k\geq 2$.

(i) If $\mu_{n-1}(G) > \frac{(k-1)n\Delta}{(n-k+1)(k-1)+4(\delta-k+2)(n-\delta-1)}$, then $\kappa(G)\geq k$.

(ii) If $G$ is triangle-free and $\mu_{n-1}(G) > \frac{(k-1)n\Delta}{(n-k+1)(k-1)+4(2\delta-k+1)(n-2\delta)}$, then $\kappa(G)\geq k$.
\end{thm}


For edge-connectivity, in this paper we obtain the following two theorems, where Theorem \ref{thm-girth-kappa'}
improves Theorems \ref{thm-LiLT13} and \ref{thm-LiLT19}, and Theorem \ref{thm-clique-kappa'}
extends Theorem \ref{thm-girth-kappa'} when $g(G)=3$.

\begin{thm}\label{thm-girth-kappa'}
Let $k$ be an integer and $G$ be a connected graph of order $n$ with minimum degree $\delta\geq k\geq 2$ and girth $g\geq 3$.
If $\mu_{n-1}(G) \geq  \frac{(k-1)n}{N(\delta,g)(n-N(\delta,g))}$, then $\kappa'(G)\geq k$.
\end{thm}

\begin{thm}\label{thm-clique-kappa'}
Let $r\geq 2$ and $k$ be integers, and $G$ be a connected graph of order $n$
with minimum degree $\delta\geq k\geq 2$ and clique number $\omega(G)\leq r$.
Let $\varphi(\delta,r) = \max\{\delta+1,\lfloor\frac{r\delta}{r-1}\rfloor\}$.
If $\mu_{n-1}(G) \geq  \frac{(k-1)n}{\varphi(\delta,r)(n-\varphi(\delta,r))}$, then $\kappa'(G)\geq k$.
\end{thm}

For vertex-connectivity, we obtain the following three theorems, where Theorem \ref{thm-girth-kappa}
improves Theorem \ref{thm-LLTW19} and extends Theorem \ref{thm-HoXL19} when $g(G)\geq 5$,
and Theorems \ref{thm-clique-kappa} and \ref{thm-clique-kappa2}
extend Theorem \ref{thm-girth-kappa} when $g(G)=3$.

\begin{thm}\label{thm-girth-kappa}
Let $g,k$ be integers and $G$ be a connected graph of order $n$ with
maximum degree $\Delta$, minimum degree $\delta\geq k\geq 2$ and girth $g\geq 3$.
If
$$\mu_{n-1}(G) > \frac{n(k-1)\Delta}{n(n-k+1)-(n-2N(\delta,g)+k-1)^2},$$
then $\kappa(G)\geq k$.
\end{thm}

\begin{thm}\label{thm-clique-kappa}
Let $r\geq 3$ and $k$ be integers, and $G$ be a connected graph of order $n$ with maximum degree $\Delta$,
minimum degree $\delta\geq k\geq 2$ and clique number $\omega(G)\leq r$.
Let $\phi(\delta,k,r)=\max\{(n-\frac{2(r-1)}{r-2}\delta+\frac{r(k-1)}{r-2})^2, (n-\frac{2r\delta}{r-1}+k-1)^2\}$.
If 
$$\mu_{n-1}(G) > \frac{n(k-1)\Delta}{n(n-k+1)- \phi(\delta,k,r)},$$
then $\kappa(G)\geq k$.
\end{thm}


\begin{thm}\label{thm-clique-kappa2}
Let $r\geq 2$ and $k\geq 2$ be integers, and $G$ be a connected graph of order $n$ with maximum degree $\Delta$,
minimum degree $\delta> (k-1)(r-1)$ and clique number $\omega(G)\leq r$.
If
$$\mu_{n-1}(G) > \frac{n(k-1)\Delta}{n(n-k+1)- (n-\frac{2r\delta}{r-1}+k-1)^2},$$
then $\kappa(G)\geq k$.
\end{thm}

Applying a result of Brouwer and Haemers \cite{BrHa12}, we get the following
two results for vertex-connectivity with respect to
$\mu_1(G)$ and $\mu_{n-1}(G)$.

\begin{thm}\label{thm-girth-kappa2}
Let $g,k$ be integers and $G$ be a connected graph of order $n$ with
minimum degree $\delta\geq k\geq 2$ and girth $g\geq 3$. If
\begin{equation*}
\frac{\mu_1(G)}{\mu_{n-1}(G)} < s + \sqrt{s^2-1}  \text{~or equivalently~} \frac{\mu_{n-1}(G)}{\mu_{1}(G)} > s - \sqrt{s^2-1},
\end{equation*}
then $\kappa(G)\geq k$, where $s=\frac{2(N(\delta,g)-k+1)(n-N(\delta,g))}{n(k-1)}+1$.
\end{thm}

\begin{thm}\label{thm-clique-kappa3}
Let $r\geq 2$ and $k\geq 2$ be integers, and $G$ be a connected graph of order $n$ with
minimum degree $\delta> (k-1)(r-1)$ and clique number $\omega(G)\leq r$.
If
\begin{equation*}
\frac{\mu_1(G)}{\mu_{n-1}(G)} < s + \sqrt{s^2-1}  \text{~or equivalently~} \frac{\mu_{n-1}(G)}{\mu_{1}(G)} > s - \sqrt{s^2-1},
\end{equation*}
then $\kappa(G)\geq k$, where $s=\frac{2(\frac{r}{r-1}\delta-k+1)(n-\frac{r}{r-1}\delta)}{n(k-1)}+1$.
\end{thm}

In Section 2, we display some preliminaries and mechanisms, including the bounds of Laplacian
eigenvalues and the scale of the remained connected components when deleting vertex subset or edge subset in $G$.
These will be applied in the proofs of the main results, to be presented in Section 3.
As corollaries, adjacency and signless Laplacian eigenvalue conditions which guarantee that
$G$ is $\kappa'(G)\geq k$ or $\kappa(G)\geq k$ are presented in the last section.

\section{Preliminaries}

In this section, we present some of the preliminaries to be used in the proof of main results.
For disjoint subsets $X$ and $Y$ of $V(G)$, let $E(X, Y)$ be the set of edges between
$X$ and $Y$. For $X\subseteq V(G)$, we use $d_G(X)$ or simply $d(X)$
to denote the number of edges between $X$ and $V(G)\setminus X$, that is
$d(X)=|E(X, V(G)\setminus X)|$. For a vertex $v\in V(G)$, we use $N_G(v)$
to denote the neighbor set of $v$ in $G$.
The following result is the famous theorem of Tur\'{a}n \cite{Tura41}.

\begin{lem}{(Tur\'{a}n \cite{Tura41})}\label{lem2.1}
Let $r\geq 1$ be an integer, and $G$ be a graph of order $n$.
If the clique number $\omega(G)\leq r$, then
$
|E(G)|\leq \left\lfloor\frac{r-1}{2r}\cdot n^2\right\rfloor.
$
\end{lem}

\begin{lem}\label{lem2.2}
Let $r\geq 2$ be an integer, and $G$ be a graph with minimum degree $\delta$
and clique number $\omega(G)\leq r$. Let $X$ be a nonempty proper subset of $V(G)$.
If $d(X) < \delta$, then $|X|\geq \max\{\delta+1,\left\lfloor\frac{r\delta}{r-1}\right\rfloor\}$.
\end{lem}

\begin{pf}
We first show that $X$ contains at least $\delta+1$ vertices.
Since each vertex in $X$ is adjacent to at most $|X|-1$ vertices of $X$, we obtain
$$
\delta |X|\leq \sum_{x\in X}d_G(x)\leq |X|(|X|-1) + d(X)\leq |X|(|X|-1) + \delta -1,
$$
and so $(|X|-1)(|X|-\delta)\geq 1$,
which means that $|X|\geq \delta+1$.

Next we show that $|X|\geq \left\lfloor\frac{r\delta}{r-1}\right\rfloor$.
By Lemma \ref{lem2.1}, we conclude that
\begin{equation}\label{e2.1}
|E(G[X])|\leq \frac{(r-1)|X|^2}{2r}.
\end{equation}
Since $\sum_{x\in X}d_G(x)=2|E(G[X])|+d(X)$, by (\ref{e2.1})
$$
|X|\delta \leq \sum_{x\in X}d_G(x)\leq 2\frac{(r-1)|X|^2}{2r} + d(X) \leq \frac{(r-1)|X|^2}{r} + \delta -1
$$
and so $|X|^2-\frac{r\delta}{r-1}|X|+\frac{r(\delta-1)}{r-1}\geq 0$.
It follows that
$$
(|X|-1)(|X|-\frac{r\delta}{r-1}+1)\geq \frac{1}{r-1}>0,
$$
which means that $|X|> \frac{r\delta}{r-1}-1$.
Therefore we arrive at
$|X|\geq \left\lfloor\frac{r\delta}{r-1}\right\rfloor$.
\end{pf}

\begin{lem}\label{lem2.3}
Let $r\geq 2$ be an integer, and $G$ be a graph with minimum degree $\delta\geq 2$ and clique number $\omega(G)\leq r$.
Let $S$ be a vertex-cut of $G$ and $X$ be the vertex set of a component of $G-S$.
\begin{itemize}
  \item[$(i)$] If $r\geq 3$ and $|S|<\delta$, then $|X|\geq \min\{\frac{r-1}{r-2}(\delta-|S|),\frac{r\delta}{r-1}-|S|\}$.
  \item[$(ii)$] If $r\geq 3$ and $\frac{\delta}{r-1}\leq |S|<\delta$, then $|X|\geq \frac{r-1}{r-2}(\delta-|S|)$.
  \item[$(iii)$] If $r\geq 2$ and $|S|<\frac{\delta}{r-1}$, then $|X|\geq \frac{r\delta}{r-1}-|S|$.
\end{itemize}
\end{lem}

\begin{pf}
(i) If $\omega(G[X])\leq r-1$, then by Lemma \ref{lem2.1}, we have
$2|E(G[X])|\leq \frac{r-2}{r-1}|X|^2$.
Since $\delta>|S|$, each vertex in $G[X]$ has degree at least $\delta-|S|$ and so
\begin{equation*}
|X|(\delta-|S|)\leq 2|E(G[X])| \leq \frac{r-2}{r-1}|X|^2.
\end{equation*}
Thus, in this case, we have
$|X|\geq \frac{r-1}{r-2}(\delta-|S|)$.

If $\omega(G[X])= r$, then there exists a complete subgraph $K_r$ in $G[X]$.
Consider the following two subcases.
If $\delta\leq r-1$, then $|X|\geq r\geq \delta+1$.
If $\delta> r-1$, then each vertex of $K_r$ has at least $\delta-r+1$ neighbors in $(X\cup S)\setminus V(K_r)$ and at most $r-1$ vertices of $K_r$ have common neighbors in $(X\cup S)\setminus V(K_r)$. This leads to
$|N(K_r)|\geq \frac{r(\delta-r+1)}{r-1}$ and so
$$
|X|+|S|\geq |V(K_r)|+|N(K_r)|\geq r+ \frac{r(\delta-r+1)}{r-1}=\frac{r\delta}{r-1},
$$
which implies $|X|\geq \frac{r\delta}{r-1}-|S|$.

By discussions above, we conclude that
\begin{enumerate}
  \item[{\bf (A)}] if $\delta\leq r-1$, then $|X|\geq \min\{\frac{r-1}{r-2}(\delta-|S|), \delta+1\}= \frac{r-1}{r-2}(\delta-|S|)$;
  \item[{\bf (B)}] if $\delta > r-1$, then $|X|\geq \min\{\frac{r-1}{r-2}(\delta-|S|),\frac{r\delta}{r-1}-|S|\}$.
\end{enumerate}

Combining (A) with (B), (i) is proved.

(ii) If $r\geq 3$ and $|S|\geq \frac{\delta}{r-1}$, then
\begin{equation*}
\frac{r\delta}{r-1}-|S|-\frac{r-1}{r-2}(\delta-|S|) = \frac{(r-1)|S|-\delta}{(r-1)(r-2)} \geq 0.
\end{equation*}
Therefore, by (i),
$|X|\geq \min\{\frac{r-1}{r-2}(\delta-|S|),\frac{r\delta}{r-1}-|S|\}= \frac{r-1}{r-2}(\delta-|S|)$.

(iii) If $r\geq 3$ and $|S|<\frac{\delta}{r-1}$, then
\begin{equation*}
\frac{r\delta}{r-1}-|S|-\frac{r-1}{r-2}(\delta-|S|) = \frac{(r-1)|S|-\delta}{(r-1)(r-2)}<0 .
\end{equation*}
Therefore, by (i), $|X|\geq \min\{\frac{r-1}{r-2}(\delta-|S|),\frac{r\delta}{r-1}-|S|\}=\frac{r\delta}{r-1}-|S|$.

If $r=2$ and $|S|<\delta$, then $X$ contains at least two vertices and there exists one edge $xy$ in $G[X]$.
As $r=2$, $G$ is triangle-free and so $N(x)\cap N(y) =\emptyset$.
Since $N(x)\cup N(y)\subseteq X\cup S$, it follows that
$$
|X|+|S|=|X\cup S|\geq |N(x)\cup N(y)| = |N(x)| + |N(y)|\geq 2\delta
$$
and thus $|X|\geq 2\delta-|S|=\frac{r\delta}{r-1}-|S|$.
The result follows.
\end{pf}

\bigskip
For any two vertices $u,v$ in $G$, let $d(u,v)$ be the length of a shortest path
between $u$ and $v$ in $G$. For any nonempty set $S\subseteq V$, let
$d(v,S)=\min\{d(v,w), \forall w\in S\}$ for any vertex $v\in V(G)$. In particular,
if $v\in S$, then $d(v,S)=0$.

\begin{lem}\label{lem2.4}
Let $G$ be a simple connected graph with minimum degree $\delta\geq 2$ and girth $g\geq 3$.
Let $S$ be a vertex-cut of $G$ and $X$ be the vertex set of a component of $G-S$.
If $|S|<\delta$, then $|X|\geq N(\delta,g)-|S|$.
\end{lem}

\begin{pf}
{\bf Claim 1.} $X$ contains at least $\delta+1-|S|$ vertices.

\medskip
Since each vertex in $X$ is adjacent to at most $|X|-1$ vertices of $X$ and at most
$|S|$ vertices of $S$, we obtain
$$
\delta |X|\leq \sum_{x\in X}d_G(x)\leq |X|(|X|-1+|S|),
$$
and so $|X|\geq \delta+1-|S|$.
Thus Claim 1 holds and implies that $|X|\geq 2$.

\medskip
{\bf Claim 2.} There exists a vertex $v\in X$ such that $d(v,S)\geq t$.

\medskip
If $t=1$, then Claim 2 holds obviously. So we only need to consider $t\geq 2$.
Suppose to the contrary that each vertex $v\in X$ satisfies $d(v,S)\leq t-1$.
Let $v_0$ be an arbitrary vertex in $X$ and $\{v_1,v_2,\dots,v_{\delta}\}\subseteq N(v_0)$ be the subset of
the neighbors of $v_0$ in $G$. For each $i\in \{1,2,\dots,\delta\}$, let $P_i$ be a shortest path
from $v_i$ to $S$, then $|E(P_i)|\leq t-1$. Note that $v_i$ may be in $S$ and $P_i$ may be trivial.
Since $|S|\leq \delta-1$, there exist at least two paths $P_j$ and $P_k$ with
$1\leq j< k\leq \delta$ such that $V(P_j)\cap V(P_k)\neq \emptyset$. Thus,
$P_j\cup P_k\cup \{v_0v_j,v_0v_k\}$ contains a cycle $C$ of length
$$
\ell(C)\leq |E(P_j)|+|E(P_k)|+2\leq 2t<g,
$$
a contradiction to the girth of $G$ is $g$. Claim 2 is proved.

\medskip
(i) Assume that $g=2t+1$ is odd and $v\in X$ with $d(v,S)\geq t$.
Then $N_i(v)\subseteq X\cup S$ for each $0\leq i\leq t$, where $N_i(v)=\{u\in V(G): d(u,v)=i\}$.
Furthermore, for each $1\leq i\leq t-1$ and for any distinct vertices $x,y \in N_i(v)$,
the neighbors of $x$ and $y$ in $N_{i+1}(v)$ are distinct
as $G[X]$ contains no cycle of length less than $g$. Hence,
  \begin{eqnarray*}
   |X| + |S|=|X\cup S| &\geq& |N_0(v)|+|N_1(v)|+|N_2(v)|+\cdots+|N_t(v)| \\
        &\geq& 1+\delta+\delta(\delta-1)+\cdots+\delta(\delta-1)^{t-1}\\
        &=&1+\delta\sum_{i=0}^{t-1}(\delta-1)^i=N(\delta,g).
  \end{eqnarray*}

\medskip
(ii) Assume that $g=2t+2$ is even and $v\in X$ with $d(v,S)\geq t$.
Let $\{v_1,v_2,\dots,v_{\delta}\}\subseteq N(v)$ be the subset of the neighbors of $v$.
Without loss of generality, assume that $P$ is the shortest path from $v$ to $v'\in S$ passing $v_1$
and $P_1$ is the subpath of $P$ from $v_1$ to $S$.
Let $P_i$ be a shortest path from $v_i$ to $S$ for each $i\in \{2,3,\dots,\delta\}$.

\medskip
{\bf Claim 3.} There exists a neighbor $u\in X$ of $v$ such that $d(u,S)\geq t$.

\medskip
Suppose that $d(v_i,S)\leq t-1$ for each $2\leq i\leq \delta$.
If there exists some $i\geq 2$ such that $V(P_i)\cap V(P_1)\neq \emptyset$, then
$$
t-1\geq d(v_i,S)= |E(P_i)|\geq |E(P_1)|=|E(P)|-1\geq t-1
$$
and so $P_i\cup P_1\cup \{vv_1, vv_i\}$ contains a cycle $C$ of length $\ell(C)\leq 2t$, a contradiction.
In this case, if $\delta=2$, then $|S|=1$ and $V(P_2)\cap V(P_1)\neq \emptyset$, 
which yields a contradiction to $g>2t$.
Hence, Claim 3 is true for $\delta=2$. Next, it suffices to consider $\delta\geq 3$.
If $V(P_i)\cap V(P_1) = \emptyset$ for each $2\leq i\leq \delta$, then
there exist at least two paths $P_i$ and $P_j$ with
$2\leq i< j\leq \delta$ such that $V(P_i)\cap V(P_j)\neq \emptyset$ as $|S\setminus \{v'\}|\leq \delta-2$.
Thus, $P_i\cup P_j\cup \{vv_i,vv_j\}$ contains a cycle $C$ of length $\ell(C)\leq 2t$, a contradiction.
This completes the proof of Claim 3.

By Claim 3, assume that $u$ is a neighbor of $v$ such that $d(u,S)\geq t$. Then
$N_i(uv)\subseteq X\cup S$ for each $1\leq i\leq t$, where $N_i(uv)=\{w\in V\setminus \{u,v\}: d(w,\{u,v\})=i\}$.
Furthermore, for each $1\leq i\leq t-1$ and for any distinct vertices $x,y \in N_i(uv)$,
the neighbors of $x$ and $y$ in $N_{i+1}(uv)$ are distinct
and $N(u)\cap N(v)=\emptyset$ as $g(G[X\cup S])\geq g=2t+2$. Hence,
  \begin{eqnarray*}
    |X| + |S|=|X\cup S|  &\geq& 2+|N_1(uv)|+|N_2(uv)|+\cdots+|N_t(uv)| \\
        &\geq& 2+2(\delta-1)+2(\delta-1)(\delta-1)+\cdots+2(\delta-1)(\delta-1)^{t-1} \\
        &=&2\sum_{i=0}^{t}(\delta-1)^i=N(\delta,g).
  \end{eqnarray*}
The result follows.
\end{pf}


\begin{lem}\label{lem2.5}
Let $G$ be a simple connected graph with minimum degree $\delta\geq 2$ and girth $g\geq 3$,
$X$ be a non-empty proper subset of $V(G)$.
If $d(X)<\delta$, then $|X|\geq N(\delta,g)$.
\end{lem}

\begin{pf}
Let $F$ be the set of edges between $X$ and $V(G)\setminus X$, and $S$ be the set of end-vertices
of $F$ in $X$, that is $S=V(F)\cap X$. Since $d(X)<\delta$, by Lemma \ref{lem2.2} we have
$|X|\geq \delta+1$. Thus, $X\setminus S\neq \emptyset$ and so $S$ is a vertex cut of $G$ with $|S|\leq d(X)< \delta$.
Let $X_1,\dots,X_k\subseteq X$ be the vertex sets of the components of $G-S$, where $k\geq 1$.
By Lemma \ref{lem2.4}, $|X_1|\geq N(\delta,g)-|S|$ and so
$|X|\geq |X_1|+|S|\geq N(\delta,g)$.
\end{pf}

\begin{cor}\label{cor2.5'}
Let $G$ be a simple graph of order $n$ with minimum degree $\delta\geq 2$ and girth $g\geq 3$.

(i) If $n<2N(\delta,g)-\kappa(G)$, then $\kappa(G)=\delta(G)$.

(ii) If $n<2N(\delta,g)$, then $\kappa'(G)=\delta(G)$.
\end{cor}

\begin{pf}
(i) Suppose to the contrary that $\kappa(G)<\delta(G)$.
Assume that $S$ is a minimum vertex-cut of $G$ and $X$ is the
vertex set of a minimum component of $G-S$. Let $Y=V(G)-(X\cup S)$.
By Lemma \ref{lem2.4}, $|Y|\geq |X|\geq N(\delta,g)-\kappa(G)$ and so
$n=|X|+|Y|+|S|\geq 2N(\delta,g)-\kappa(G)$, which is a contradiction.

(ii) Suppose to the contrary that $\kappa'(G)<\delta(G)$.
Assume that $F=E(X,Y)$ is a minimum edge-cut of $G$ and $|Y|\geq |X|$.
By Lemma \ref{lem2.5}, $|Y|\geq |X|\geq N(\delta,g)$ and so
$n=|X|+|Y|\geq 2N(\delta,g)$, which is a contradiction.
\end{pf}

\bigskip
Let $x=(x_1,x_2,\dots,x_n)^T\in \mathbb{R}^n$, and let $G$ be a
graph with vertex set $V(G)=\{1,2,\dots, n\}$. Then $x$ can be considered
as a function defined on $V(G)$, that is, for any
vertex $i$, we map it to $x_i=x(i)$. Fiedler \cite{Fied75} derived a very useful
expression for algebraic connectivity as follows.

\begin{lem}{(Fiedler \cite{Fied75})}\label{lem-Fied75}
Let $G$ be a connected graph with vertex set $V=\{1,2,\dots,n\}$ and edge set $E$. Then
the algebraic connectivity of $G$ is positive and
\begin{equation*}
\mu_{n-1}(G) =  \min_x f(x)= \min_x\frac{n\sum\limits_{ij\in E}(x_i-x_j)^2}{\sum\limits_{i,j\in V, i<j}(x_i-x_j)^2},
\end{equation*}
where the minimum is taken over all non-constant vectors $x=(x_1,x_2,\dots,x_n)^T\in \mathbb{R}^n$.
The characteristic vectors $y=(y_1,y_2,\dots,y_n)^T$ of $G$ corresponding to $\mu_{n-1}(G)$ are
then those non-constant vectors for which the minimum of $f(x)$ is attained and for which $\sum_{i=1}^n y_i=0$.
\end{lem}

\begin{lem}\label{lem2.7}
Let $G=(V,E)$ be a graph, and $X$ be a nonempty proper subset of $V$ and $Y=V\setminus X$. Then
\begin{equation*}
\mu_{n-1}(G) \leq \frac{nd(X)}{|X||Y|}.
\end{equation*}
\end{lem}

\begin{pf}
Let $x=(x_1,x_2,\dots,x_n)^T$ be a real vector. If $i\in X$, then set $x_i=1$; if $i\in Y$, then set $x_i=-1$. By Lemma \ref{lem-Fied75},
\begin{equation}\label{ee2.4}
\mu_{n-1}(G)\leq \frac{n\sum\limits_{ij\in E}(x_i-x_j)^2}{\sum\limits_{i,j\in V, i<j}(x_i-x_j)^2}
\end{equation}
holds for the real vector $x$.
Applying the values of the entries of $x$ into the inequality (\ref{ee2.4}), we obtain
\begin{equation*}
\sum\limits_{ij\in E}(x_i-x_j)^2 = \sum\limits_{ij\in E(X,Y)}(1-(-1))^2 = 4d(X),
\end{equation*}
\begin{equation*}
\sum\limits_{i,j\in V, i<j}(x_i-x_j)^2  = \sum\limits_{i\in X, j\in Y}(1-(-1))^2=4|X||Y|.
\end{equation*}
By (\ref{ee2.4}), the result follows.
\end{pf}

\begin{lem}\label{lem2.8}
Let $G=(V,E)$ be a graph of order $n$, and $S$ be an arbitrary minimum vertex-cut of $G$
and $X$ be the vertex set of a component of $G-S$, and $Y=V-(S\cup X)$. Then
\begin{equation*}
\mu_{n-1}(G) \leq \frac{n d(S)}{n(n-|S|)-(|X|-|Y|)^2}.
\end{equation*}
\end{lem}

\begin{pf}
Let $x=(x_1,x_2,\dots,x_n)^T$ be a real vector. If $i\in X$, then set $x_i=1$; if $i\in Y$, then set $x_i=-1$;
if $i\in S$, then set $x_i=0$. By Lemma \ref{lem-Fied75},
\begin{equation}\label{ee2.7}
\mu_{n-1}(G)\leq \frac{n\sum\limits_{ij\in E}(x_i-x_j)^2}{\sum\limits_{i,j\in V, i<j}(x_i-x_j)^2}
\end{equation}
holds for the real vector $x$.
Applying the values of the entries of $x$ into the inequality (\ref{ee2.7}), we obtain
\begin{equation}\label{ee2.8}
\sum\limits_{ij\in E}(x_i-x_j)^2 = \sum\limits_{ij\in E(S,X\cup Y)}(x_i-x_j)^2 = \sum\limits_{ij\in E(S,X\cup Y)}1 = d(S),
\end{equation}
\begin{align}\label{ee2.9}
\sum\limits_{i,j\in V, i<j}(x_i-x_j)^2
& = \sum\limits_{i\in X, j\in S}(x_i-x_j)^2 + \sum\limits_{i\in Y, j\in S}(x_i-x_j)^2 + \sum\limits_{i\in X, j\in Y}(x_i-x_j)^2 \nonumber\\
& = \sum\limits_{i\in X, j\in S}(1-0)^2 + \sum\limits_{i\in Y, j\in S}((-1)-0)^2 + \sum\limits_{i\in X, j\in Y}(1-(-1))^2 \nonumber\\
& = |S||X|+|S||Y|+4|X||Y|  \nonumber\\
& = (n-|X|-|Y|)(|X|+|Y|)+4|X||Y|  \nonumber\\
& = n(n-|S|) - (|X|-|Y|)^2.
\end{align}
Substituting (\ref{ee2.8}) and (\ref{ee2.9}) in (\ref{ee2.7}), the result follows.
\end{pf}

\begin{lem}{(Haemers \cite{Haem95})}\label{lem-Haem95}
Let $G$ be a graph on $n$ vertices, and let $X$ and $Y$ be
disjoint sets of vertices, such that there is no edge between $X$ and $Y$. Then
$$
\frac{|X||Y|}{(n-|X|)(n-|Y|)}\leq \left(\frac{\mu_1(G)-\mu_{n-1}(G)}{\mu_1(G)+\mu_{n-1}(G)}\right)^2.
$$
\end{lem}

For applications, a useful Lemma can be derived from Lemma \ref{lem-Haem95} as follows.

\begin{lem}{(Brouwer and Haemers \cite{BrHa12})}\label{lem-BrHa12}
Let $G$ be a connected graph on $n$ vertices, and let $X$ and $Y$ be
disjoint sets of vertices, such that there is no edge between $X$ and $Y$. Then
$$
\frac{|X||Y|}{n(n-|X|-|Y|)}\leq \frac{(\mu_1(G)-\mu_{n-1}(G))^2}{4\mu_1(G)\mu_{n-1}(G)}.
$$
\end{lem}

\section{The proof of main results}

\noindent\textbf{\bf{Proof of Theorem \ref{thm-girth-kappa'}.}}
To the contrary, suppose that $1\leq \kappa'(G)\leq k-1$.
Let $F$ be an arbitrary minimum edge-cut of $G$,
and $X$, $Y$ be the vertex sets of two components of $G-F$ with $|X|\leq|Y|$. Thus $d(X)=\kappa'(G)\leq k-1$.
By Lemma \ref{lem2.5} and $d(X)<\delta$, we obtain $|X| \geq N(\delta,g)$.
Since $|Y|\geq |X|$ and $|X|+|Y|=n$,
\begin{equation}\label{e3.1}
|X|\cdot|Y|\geq N(\delta,g)(n-N(\delta,g)).
\end{equation}
By Lemma \ref{lem2.7} and (\ref{e3.1}), we have
$$
\mu_{n-1}(G)\leq \frac{nd(X)}{|X||Y|} \leq \frac{(k-1)n}{N(\delta,g)(n-N(\delta,g))}.
$$
According to the hypothesis, it follows that
$\mu_{n-1}(G)=\frac{nd(X)}{|X||Y|}=\frac{(k-1)n}{N(\delta,g)(n-N(\delta,g))}$.
By the proof of Lemma \ref{lem2.7},
$\mu_{n-1}(G) = \frac{n\sum\nolimits_{ij\in E}(x_i-x_j)^2}{\sum\nolimits_{i,j\in V, i<j}(x_i-x_j)^2}$,
where $x_i=1$ if $i\in X$ and $x_i=-1$ if $i\in Y$.
By Lemma \ref{lem-Fied75}, $x$ is a characteristic vector of $G$ corresponding to $\mu_{n-1}(G)$.
Since $d(X)<\delta$ and $|X|\geq N(\delta,g)\geq \delta+1$, there exists one vertex $j$ in $X$
such that its neighbor set $N_G(j)\subset X$. Thus, by $\mu_{n-1}(G)x=(D-A)x$, we have
$\mu_{n-1}(G)x_j=|N_G(j)|x_j-\sum_{\ell\in N_G(j)}x_{\ell}$.
Since $x_j=x_\ell=1$, it indicates $\mu_{n-1}(G)=0$ and so $k-1=0$,
which is a contradiction to $k\geq 2$. Hence, $\kappa'(G)\geq k$.
\hfill\rule{2mm}{2mm}

\begin{rem}{\rm
The result in Theorem \ref{thm-girth-kappa'} improves the one of Theorem \ref{thm-LiLT13}
when $\delta\geq 3$ and improves the one of Theorem \ref{thm-LiLT19} when $\delta\geq 3$ and $g\geq 5$.
In fact, if $n< 2N(\delta,g)$, then by Corollary \ref{cor2.5'} we have $\kappa'(G)=\delta(G)$.
Therefore, we only need to compare the bounds when $n\geq 2N(\delta,g)$. Note that
$N(\delta,g)>N(\delta,g)-\sum_{i=1}^{t-1}(\delta-1)^i=f(\delta,g)$ when $\delta\geq 3$ and
$g\geq 5$, and $N(\delta,g)>\frac{4}{9}N(\delta,g)$. As $N(\delta,g)\leq \frac{n}{2}$, it follows that
$N(\delta,g)(n-N(\delta,g))>\frac{4}{9}N(\delta,g)(n-\frac{4}{9}N(\delta,g))$ and
$N(\delta,g)(n-N(\delta,g))>f(\delta,g)(n-f(\delta,g))$, and so
$$
\frac{(k-1)n}{N(\delta,g)(n-N(\delta,g))}<\frac{(k-1)n}{\frac{4}{9}N(\delta,g)(n-\frac{4}{9}N(\delta,g))}, \
\frac{(k-1)n}{N(\delta,g)(n-N(\delta,g))}<\frac{(k-1)n}{f(\delta,g)(n-f(\delta,g))}.
$$
}
\end{rem}

\noindent\textbf{\bf{Proof of Theorem \ref{thm-clique-kappa'}.}}
To the contrary, suppose that $1\leq \kappa'(G)\leq k-1$.
Let $F$ be an arbitrary minimum edge-cut of $G$,
and $X$, $Y$ be the vertex sets of two components of $G-F$ with $|X|\leq|Y|$. Thus $d(X)=\kappa'(G)\leq k-1$.
By Lemma \ref{lem2.2} and $d(X)<\delta$, we obtain
$|X| \geq \varphi(\delta,r)=\max\{\delta+1,\lfloor\frac{r\delta}{r-1}\rfloor\}$.
Since $|Y|\geq |X|$ and $|X|+|Y|=n$,
\begin{equation}\label{e3.2}
|X|\cdot|Y|\geq \varphi(\delta,r)(n-\varphi(\delta,r)).
\end{equation}
By Lemma \ref{lem2.7} and (\ref{e3.2}), we have
$$
\mu_{n-1}(G)\leq \frac{nd(X)}{|X||Y|} \leq \frac{(k-1)n}{\varphi(\delta,r)(n-\varphi(\delta,r))}.
$$
According to the hypothesis, it follows that
$\mu_{n-1}(G)=\frac{nd(X)}{|X||Y|}=\frac{(k-1)n}{\varphi(\delta,r)(n-\varphi(\delta,r))}$.
By the proof of Lemma \ref{lem2.7},
$\mu_{n-1}(G) = \frac{n\sum\nolimits_{ij\in E}(x_i-x_j)^2}{\sum\nolimits_{i,j\in V, i<j}(x_i-x_j)^2}$,
where $x_i=1$ if $i\in X$ and $x_i=-1$ if $i\in Y$.
By Lemma \ref{lem-Fied75}, $x$ is a characteristic vector of $G$ corresponding to $\mu_{n-1}(G)$.
Since $d(X)<\delta$ and $|X|\geq \varphi(\delta,r)\geq \delta+1$, there exists one vertex $j$ in $X$
such that $N_G(j)\subset X$. Thus, by $\mu_{n-1}(G)x=(D-A)x$, we have
$\mu_{n-1}(G)x_j=|N_G(j)|x_j-\sum_{\ell\in N_G(j)}x_{\ell}$.
Since $x_j=x_\ell=1$, it indicates $\mu_{n-1}(G)=0$ and so $k-1=0$,
which is a contradiction to $k\geq 2$. Hence, $\kappa'(G)\geq k$.
\hfill\rule{2mm}{2mm}

\bigskip
\noindent\textbf{\bf{Proof of Theorem \ref{thm-girth-kappa}.}}
To the contrary, suppose that $1\leq \kappa=\kappa(G)\leq k-1$.
Let $S$ be an arbitrary minimum vertex-cut
and $X$ be the vertex set of a minimum component of $G-S$, and $Y=V-(S\cup X)$.
By Lemma \ref{lem2.4} and $|S|=\kappa\leq k-1<\delta$, we obtain $|X|\geq N(\delta,g)-|S|$. Thus
\begin{equation}\label{e3.3}
N(\delta,g)-k+1 \leq |X|\leq |Y| \leq n-N(\delta,g),
\end{equation}
and so $(|X|-|Y|)^2\leq (n-2N(\delta,g)+k-1)^2$. Therefore,
\begin{equation}\label{e3.5}
n(n-|S|)-(|X|-|Y|)^2  \geq  n(n-k+1)-(n-2N(\delta,g)+k-1)^2.
\end{equation}
By $N(\delta,g)\geq \delta+1>k$ and (\ref{e3.3}), we have $n-k+1>n-2N(\delta,g)+k-1\geq 0$, which implies
$n(n-k+1)-(n-2N(\delta,g)+k-1)^2>0$.
Combining Lemma \ref{lem2.8} with inequality (\ref{e3.5}), we have
$$
\mu_{n-1}(G)\leq \frac{n d(S)}{n(n-|S|)-(|X|-|Y|)^2}  \leq \frac{n(k-1)\Delta}{n(n-k+1)-(n-2N(\delta,g)+k-1)^2},
$$
which is a contradiction to the hypothesis. Hence, $\kappa(G)\geq k$.
\hfill\rule{2mm}{2mm}

\begin{rem}{\rm
The result in Theorem \ref{thm-girth-kappa} improves the one of Theorem \ref{thm-LLTW19}
and extends the one of Theorem \ref{thm-HoXL19} when $g\geq 5$.
In fact, if $n< 2N(\delta,g)-\kappa(G)$, then by Corollary \ref{cor2.5'} we have $\kappa(G)=\delta(G)$.
Therefore, we only need to compare the bounds when $n\geq 2N(\delta,g)-\kappa(G)$.

(i) Theorem \ref{thm-girth-kappa} improves Theorem \ref{thm-LLTW19}.
Denote $N:=N(\delta,g)$, $\kappa:=\kappa(G)$ and $\nu:=\nu(\delta,g,k)$.
Then $n\geq 2N-k+1>N$ and so $n-k+1\geq 2(N-k+1)$.
As $\nu=N-(k-1)\sum\nolimits_{i=0}^{t-1}(\delta-1)^i\leq N-k+1$,
we get $n>2(N-k+1)\geq 2\nu$. Hence,
\begin{align*}
n(n-k+1)-(n-2N+k-1)^2
& = (n-k+1)(k-1)+4(N-k+1)(n-N)\\
& \geq  2(N-k+1)(k-1)+4(N-k+1)(n-N)\\
& = 2(N-k+1)(n-(N-k+1)+(n-N))\\
& > 2(N-k+1)(n-(N-k+1))\geq 2\nu(n-\nu),
\end{align*}
and we arrive at
$\frac{n(k-1)\Delta}{n(n-k+1)-(n-2N(\delta,g)+k-1)^2} < \frac{n(k-1)\Delta}{2\nu(n-\nu)}.$

(ii) Theorem \ref{thm-girth-kappa} extends Theorem \ref{thm-HoXL19}.
Suppose $n\geq 2N(\delta,g)-k+1$ and $\delta\geq 2$.
If $g\geq 3$, then $N(\delta,g)\geq N(\delta,3)=\delta+1$ and so
\begin{align*}
n(n-k+1)-(n-2N(\delta,g)+k-1)^2
&\geq n(n-k+1)-(n-2(\delta+1)+k-1)^2 \\
&= (n-k+1)(k-1)+4(\delta-k+2)(n-\delta-1).
\end{align*}
If $G$ is triangle-free, then
$g\geq 4$ and $N(\delta,g)\geq N(\delta,4)=2\delta$, and thus
\begin{align*}
n(n-k+1)-(n-2N(\delta,g)+k-1)^2
&\geq n(n-k+1)-(n-4\delta +k-1)^2 \\
&= (n-k+1)(k-1)+4(2\delta-k+1)(n-2\delta).
\end{align*}
Therefore, the lower bound on $\mu_{n-1}(G)$ in Theorem \ref{thm-girth-kappa}
is less than or equal to the one in Theorem \ref{thm-HoXL19}, and Theorem \ref{thm-girth-kappa}
extends Theorem \ref{thm-HoXL19} when $g(G)\geq 5$.
}
\end{rem}

\bigskip
\noindent\textbf{\bf{Proof of Theorem \ref{thm-clique-kappa}.}}
To the contrary, suppose that $1\leq \kappa=\kappa(G)\leq k-1$.
Let $S$ be an arbitrary minimum vertex-cut
and $X$ be the vertex set of a minimum component of $G-S$, and $Y=V-(S\cup X)$.
Consider the following two cases.

(i) $\frac{\delta}{r-1}\leq |S|=\kappa<\delta$. By Lemma \ref{lem2.3} (ii),
\begin{equation}\label{e3.6}
\frac{r-1}{r-2}(\delta-\kappa)\leq |X|\leq |Y| \leq n-\kappa-\frac{r-1}{r-2}(\delta-\kappa),
\end{equation}
and so $n-\frac{2(r-1)}{r-2}\delta+\frac{r\kappa}{r-2}\geq 0$ and $(|X|-|Y|)^2\leq (n-\frac{2(r-1)}{r-2}\delta+\frac{r\kappa}{r-2})^2$. Therefore,
\begin{align}\label{e3.7}
n(n-|S|)-(|X|-|Y|)^2& \geq n(n-\kappa)-(n-\frac{2(r-1)}{r-2}\delta+\frac{r\kappa}{r-2})^2 \nonumber\\
& \geq n(n-k+1)-(n-\frac{2(r-1)}{r-2}\delta+\frac{r(k-1)}{r-2})^2.
\end{align}
By $\delta>k-1$ and (\ref{e3.6}), we have $n-k+1>n-\frac{2(r-1)}{r-2}\delta+\frac{r(k-1)}{r-2}\geq 0$, which implies
$n(n-k+1)-(n-\frac{2(r-1)}{r-2}\delta+\frac{r(k-1)}{r-2})^2>0$.
Combining (\ref{e3.7}) with $d(S)\leq (k-1)\Delta$, by Lemma \ref{lem2.8}, we have
\begin{align}\label{e3.8}
\mu_{n-1}(G)\leq \frac{n(k-1)\Delta}{n(n-k+1)- (n-\frac{2(r-1)}{r-2}\delta+\frac{r(k-1)}{r-2})^2}.
\end{align}

(ii) $|S|=\kappa<\frac{\delta}{r-1}$.
By Lemma \ref{lem2.3} (iii), we get
\begin{equation}\label{ee3.6}
\frac{r\delta}{r-1}-k+1\leq \frac{r\delta}{r-1}-\kappa\leq |X|\leq |Y| \leq n-\frac{r\delta}{r-1},
\end{equation}
and so $(|X|-|Y|)^2\leq (n-\frac{2r\delta}{r-2}+k-1)^2$. Therefore,
\begin{equation}\label{ee3.7}
n(n-|S|)-(|X|-|Y|)^2 \geq n(n-k+1)-(n-\frac{2r\delta}{r-1}+k-1)^2.
\end{equation}
By $\delta>k-1$ and (\ref{ee3.6}), we have $n-k+1>n-\frac{2r\delta}{r-1}+k-1\geq 0$, which implies
$n(n-k+1)-(n-\frac{2r\delta}{r-1}+k-1)^2>0$.
Combining (\ref{ee3.7}) with $d(S)\leq (k-1)\Delta$, by Lemma \ref{lem2.8}, we have
\begin{align}\label{ee3.8}
\mu_{n-1}(G)\leq \frac{n(k-1)\Delta}{n(n-k+1)- (n-\frac{2r\delta}{r-1}+k-1)^2}.
\end{align}

Now, let $\phi(\delta,k,r)=\max\{(n-\frac{2(r-1)}{r-2}\delta+\frac{r(k-1)}{r-2})^2, (n-\frac{2r\delta}{r-1}+k-1)^2\}$.
By (\ref{e3.8}) and (\ref{ee3.8}), we have
$$
\mu_{n-1}(G) \leq\frac{n(k-1)\Delta}{n(n-k+1)-\phi(\delta,k,r)},
$$
which is a contradiction to the hypothesis. Hence, $\kappa(G)\geq k$.
\hfill\rule{2mm}{2mm}

\bigskip
\noindent\textbf{\bf{Proof of Theorem \ref{thm-clique-kappa2}.}}
If $r=2$, then $g(G)\geq 4$ and so $N(\delta,g)\geq N(\delta,4)=2\delta=\frac{r}{r-1}\delta$.
Thus, by Theorem \ref{thm-girth-kappa}, the theorem holds when $r=2$.
Next we consider $r\geq 3$.
To the contrary, suppose that $1\leq \kappa=\kappa(G)\leq k-1$.
Let $S$ be an arbitrary minimum vertex-cut
and $X$ be the vertex set of a minimum component of $G-S$, and $Y=V-(S\cup X)$.
By Lemma \ref{lem2.3} and $|S|=\kappa\leq k-1<\frac{\delta}{r-1}$, we obtain
$|X|\geq \frac{r\delta}{r-1}-\kappa$. Thus,
$$
\frac{r\delta}{r-1}-k+1\leq \frac{r\delta}{r-1}-\kappa\leq |X|\leq |Y| \leq n-\frac{r\delta}{r-1},
$$
Using a similar argument as in the proof of Theorem \ref{thm-clique-kappa}, we have
$$
\mu_{n-1}(G)\leq \frac{n(k-1)\Delta}{n(n-k+1)- (n-\frac{2r\delta}{r-1}+k-1)^2}.
$$
which is a contradiction to the hypothesis. Hence, $\kappa(G)\geq k$ and the result follows.
\hfill\rule{2mm}{2mm}

\begin{rem}{\rm
If $\omega(G)\geq 3$, then $g(G)=3$. In this case, since $\frac{r\delta}{r-1}\geq \delta+1$ for $2\leq r\leq \delta+1$
and $\frac{r\delta}{r-1}> \delta+1$ for $3\leq r\leq \delta$, it follows that
Theorem \ref{thm-clique-kappa'} extends Theorem \ref{thm-girth-kappa'}, and
Theorems \ref{thm-clique-kappa} and \ref{thm-clique-kappa2} extends Theorem \ref{thm-girth-kappa} when $g(G)=3$.
}
\end{rem}

\bigskip
\noindent\textbf{\bf{Proof of Theorem \ref{thm-girth-kappa2}.}}
To the contrary, suppose that $1\leq \kappa=\kappa(G)\leq k-1$.
Let $S$ be an arbitrary minimum vertex-cut
and $X$ be the vertex set of a minimum component of $G-S$, and $Y=V-(S\cup X)$.
By Lemma \ref{lem2.4} and $1\leq \kappa\leq k-1<\delta$, we obtain
\begin{equation*}
N(\delta,g)-\kappa\leq |X|\leq |Y|\leq n-N(\delta,g),
\end{equation*}
and so
\begin{equation*}
|X|\cdot |Y|\geq (N(\delta,g)-\kappa)(n-N(\delta,g))\geq (N(\delta,g)-k+1)(n-N(\delta,g)).
\end{equation*}
Combining this with $n-|X|-|Y|=\kappa \leq k-1$, by Lemma \ref{lem-BrHa12},
\begin{equation}\label{e3.11}
\frac{(\mu_1(G)-\mu_{n-1}(G))^2}{4\mu_1(G)\mu_{n-1}(G)} \geq \frac{|X||Y|}{n(n-|X|-|Y|)}\geq \frac{(N(\delta,g)-k+1)(n-N(\delta,g))}{n(k-1)}.
\end{equation}
Set $t=\frac{\mu_1(G)}{\mu_{n-1}(G)}$ and $s=\frac{2(N(\delta,g)-k+1)(n-N(\delta,g))}{n(k-1)}+1$.
Substituting $t$ and $s$ in (\ref{e3.11}), we obtain
$t+t^{-1}\geq 2s.$
Since $t\geq 1$ and $s\geq 1$, $t \geq s + \sqrt{s^2-1}$ is necessary.
This contradicts to the hypothesis. Therefore, $\kappa(G)\geq k$.
\hfill\rule{2mm}{2mm}

\bigskip
\noindent\textbf{\bf{Proof of Theorem \ref{thm-clique-kappa3}.}}
To the contrary, suppose that $1\leq \kappa=\kappa(G)\leq k-1$.
Let $S$ be an arbitrary minimum vertex-cut
and $X$ be the vertex set of a minimum component of $G-S$, and $Y=V-(S\cup X)$.
By Lemma \ref{lem2.3} (iii) and $|S|\leq k-1<\frac{\delta}{r-1}$, we obtain
$|X|\geq \frac{r\delta}{r-1}-|S|$. Thus,
$$
\frac{r\delta}{r-1}-\kappa\leq |X|\leq |Y| \leq n-\frac{r\delta}{r-1},
$$
and so
$$
|X|\cdot |Y|\geq (\frac{r\delta}{r-1}-\kappa)(n-\frac{r\delta}{r-1})\geq (\frac{r\delta}{r-1}-k+1)(n-\frac{r\delta}{r-1}).
$$
Combining this with $n-|X|-|Y|=\kappa \leq k-1$, by Lemma \ref{lem-BrHa12},
\begin{equation}\label{e3.12}
\frac{(\mu_1(G)-\mu_{n-1}(G))^2}{4\mu_1(G)\mu_{n-1}(G)} \geq \frac{|X||Y|}{n(n-|X|-|Y|)}\geq \frac{(\frac{r\delta}{r-1}-k+1)(n-\frac{r\delta}{r-1})}{n(k-1)}.
\end{equation}
Set $t=\frac{\mu_1(G)}{\mu_{n-1}(G)}$ and $s=\frac{2(\frac{r\delta}{r-1}-k+1)(n-\frac{r\delta}{r-1})}{n(k-1)}+1$.
Substituting $t$ and $s$ in (\ref{e3.12}), we obtain
$t+t^{-1}\geq 2s.$
Since $t\geq 1$ and $s\geq 1$, $t \geq s + \sqrt{s^2-1}$ is necessary.
This contradicts to the hypothesis. Therefore, $\kappa(G)\geq k$.
\hfill\rule{2mm}{2mm}

\section{Connectivity and adjacency or signless Laplacian eigenvalues}

In this section, we present the relationship between (edge-)connectivity
and the second largest adjacency eigenvalue or the second largest signless Laplacian eigenvalue.

\begin{thm}{\rm (Courant-Weyl Inequalities)}\label{thm-weyl}
Let $A$ and $B$ be Hermitian matrices of order $n$, and let $1\leq i, j\leq n$.
If $i+j\leq n+1$, then $\lambda_i(A)+\lambda_j(B) \geq \lambda_{i+j-1}(A+B)$.
\end{thm}

For real numbers $a, b$ with $b>0$ and $a\geq -b$,
let $\lambda_i(G, a, b)$ be the $i$th largest eigenvalue of the matrix
$aD + bA$.

\begin{cor}\label{cor-eigen}
Let $p\geq 0$, $b>0$ and $a\geq -b$ be real numbers and $G$ be a graph of order $n$ with minimum degree $\delta$.\\
(i) If $\lambda_2(G,a,b) < (a+b)\delta - bp$, then $\mu_{n-1}(G) > p$.  In particular,
if $q_2(G) < 2\delta - p$ or $\lambda_2(G) < \delta - p$, then $\mu_{n-1}(G) > p$.\\
(ii) If $\lambda_2(G,a,b) \leq (a+b)\delta - bp$, then $\mu_{n-1}(G) \geq p$.  In particular,
if $q_2(G) \leq 2\delta - p$ or $\lambda_2(G) \leq \delta - p$, then $\mu_{n-1}(G) \geq p$.
\end{cor}

\begin{pf}
Let $A$ and $D$ be the adjacency matrix and degree diagonal matrix of $G$.
Since $b(D-A)+(aD+bA)=(a+b)D$, by Theorem \ref{thm-weyl}, $\lambda_{n-1}(b(D-A)) + \lambda_2(aD+bA) \geq \lambda_n((a+b)D)$.
As $b>0$ and $a + b\geq 0$, $b\mu_{n-1}(G)+\lambda_2(G,a,b)\geq (a+b)\delta$. Therefore, if $\lambda_2(G,a,b) < (a+b)\delta - bp$,
then $\mu_{n-1}(G) > p$. In particular, $\lambda_2(G,1,1)=q_2(G)$ and $\lambda_2(G,0,1)=\lambda_2(G)$. Thus, (i) is proved and (ii)
can be proved similarly.
\end{pf}

\bigskip
By Corollary \ref{cor-eigen}, from the sufficient conditions on $\mu_{n-1}(G)$
in Theorems \ref{thm-girth-kappa'}-\ref{thm-clique-kappa2}, we can obtain
sufficient conditions on $\lambda_2(G,a,b)$, especially on $\lambda_2(G)$
and $q_2(G)$. For example, by Corollary \ref{cor-eigen} and Theorem \ref{thm-girth-kappa'},
we have the following corollary. Other corollaries could be stated similarly.

\begin{cor}
Let $k$ be an integer and $G$ be a connected graph of order $n$ with minimum degree $\delta\geq k\geq 2$ and girth $g\geq 3$.
If $\lambda_2(G) \leq \delta - \frac{(k-1)n}{N(\delta,g)(n-N(\delta,g))}$ or
$q_2(G) \leq 2\delta-\frac{(k-1)n}{N(\delta,g)(n-N(\delta,g))}$, then $\kappa'(G)\geq k$.
\end{cor}

\section*{Acknowledgement}

The research of Zhen-Mu Hong is supported by NNSFC (No. 11601002),
Outstanding Young Talents International Visiting Program of Anhui Provincial Department of Education (No. gxgwfx2018031)
and Key Projects in Natural Science Research of Anhui Provincial Department of Education (No. KJ2016A003).
The research of Hong-Jian Lai is supported by NNSFC (Nos. 11771039 and 11771443).

\end{document}